\documentclass{amsart}
\usepackage{epsfig}

\begin{document}

\newtheorem{thm}{Theorem}[section]
\newtheorem{lem}[thm]{Lemma}
\newtheorem{cor}[thm]{Corollary}
\newtheorem{conj}[thm]{Conjecture}

\theoremstyle{definition}
\newtheorem{defn}{Definition}[section]

\theoremstyle{remark}
\newtheorem{rmk}{Remark}[section]

\def\square{\hfill${\vcenter{\vbox{\hrule height.4pt \hbox{\vrule width.4pt
height7pt \kern7pt \vrule width.4pt} \hrule height.4pt}}}$}

\def\R{\Bbb R}
\def\Z{\Bbb Z}

\newenvironment{pf}{{\it Proof:}\quad}{\square \vskip 12pt}

\title{A degree one Borsuk-Ulam theorem}

\author{Danny Calegari}
\address{Department of Mathematics \\ UC Berkeley \\ Berkeley, CA 94704}
\email{dannyc@math.berkeley.edu}

\maketitle

\begin{abstract}
We generalize the Borsuk-Ulam theorem for maps $M^n \to \R^n$.
\end{abstract}

Everyone knows the Borsuk-Ulam theorem as a simple application of
some of the first ideas one encounters in algebraic topology. 

\begin{thm}[Borsuk-Ulam]
Let $f:S^n \to \R^n$ be any continuous map. Then there are antipodal
points in $S^n$ which are mapped to the same point under $f$.
\end{thm}

The purpose of this brief note is to observe that there is an easy
generalization of this theorem for maps $f:M^n \to \R^n$ where $M^n$ is
a closed $n$-manifold.

\begin{thm}
Let $M$ be a closed $n$-manifold. Let $f:M \to \R^n$ be any continuous
map and $g:M \to S^n$ a degree one map. Then there are points $p,q \in M$
such that $f(p) = f(q)$ and $g(p) = - g(q)$.
\end{thm}
\begin{pf}
We wiggle $g$ to be smooth and generic. By compactness of the space of
antipodal points in $S^n$, it suffices to prove the theorem in this case,
since then we can extract a subsequence of pairs of points in $M$ with
the desired properties for a sequence of degree one smooth maps
$g_i:M \to S^n$ approximating $g$.

We define the following spaces
$$\hat M \subset M \times M - \Delta = \lbrace (p,q): g(p) = -g(q) \rbrace $$
$$ S \subset S^n \times S^n - \Delta = \lbrace (p,q): p = -q \rbrace $$
Observe that $S$ is homeomorphic to $S^n$. There is an induced map
$\hat g:\hat M \to S$ given by $\hat g:(p,q) \to (g(p),g(q))$. Since
$g$ was degree one, one easily observes that there are an odd number of
points in the generic fiber of $\hat g$ so that there is some connected
component of $\hat M$ for which the restricted map $\hat g$ has odd
degree. Moreover, the $\Z/2\Z$ action on $\hat M$ and $S$ given by
interchanging the co-ordinates commutes with $\hat g$, so there is an
induced map on the quotients. We define $N = \hat M / \sim$ and call
the quotient map $h:N \to \R P^n$.

Assume on the contrary that points in $M$ mapping to antipodal points in
$S^n$ map to distinct points in $\R^n$. Then there is a map
$$\hat f:\hat M \to S^{n-1}$$
defined by
$$\hat f:(p,q) \to \frac {f(p) - f(q)} {|| f(p) - f(q) ||}$$
It is obvious that this descends to a map $j: N \to \R P^{n-1}$ where
$\R P^{n-1}$ is obtained from $S^n$ by quotienting out by the antipodal map.

In the sequel, we consider homology and cohomology with $\Z/2\Z$ coefficients.
For simplicity of notation, we omit the coefficients.

Since the degree of $h$ is odd, $h^*$ pulls back the generator $[\R P^n]$ of
$H^n(\R P^n)$ to the generator $[N]$ of $H^n(N)$. Furthermore,
if $\alpha$ generates $H^1(\R P^n)$ then
$h^*\alpha \in H^1(N)$ is an element whose $n$th power is
$[N]$. Moreover by construction for every cycle $C \in H_1(N)$
we have $h_*C \ne 0$ in $H_1(\R P^n)$ iff $j_*C \ne 0$ in 
$H_1(\R P^{n-1})$, since these are exactly the $C$ which do not lift 
to $\hat M$.

It follows that if $\beta$ denotes the generator of $H^1(\R P^{n-1})$
then $j^*\beta(C) = h^*\alpha(C)$ for all $C$, and therefore
$j^*\beta = h^*\alpha$ so that the $n$th power of
$j^*\beta$ is nontrivial. But $(j^*\beta)^n = j^*(\beta^n)$ which is
trivial, giving us a contradiction.
\end{pf}

\begin{rmk}
Notice that the proof works in exactly the same way if $g:M \to S^n$ is
a map of odd degree.
\end{rmk}

The following corollary led the author to observe the theorem above:

\begin{cor}
Let $M^n \subset \R^{n+1}$ be an embedded submanifold bounding
a closed region which contains a ball of diameter $t$. Let
$f:M^n \to \R^n$ be a continuous map. Then there are points in
$M$ at distance at least $t$ apart from each other which have the
same image under $f$.
\end{cor}
\begin{pf}
Let $g$ be the map which is radial projection of $M$ onto the boundary
of the ball of diameter $t$.
\end{pf}

\end{document}